\input amstex
\documentstyle{amsppt}
\headline={\hss\tenrm\hss\folio}
\vsize=9.0truein
\hsize=6.5truein
\parindent=20pt
\baselineskip=22pt
\NoBlackBoxes

%Sept08Version 

  \topmatter
   %TITLE
\title Dirac Operators on Non--Compact Orbifolds \endtitle 
\author Carla Farsi, \
Department of Mathematics, \ 
University of Colorado, \
395 UCB, \
Boulder, CO 80309--0395, USA. \ 
e-mail: farsi\@euclid.colorado.edu 
\endauthor
% ABSTRACT
\abstract 
 In this paper we prove that Dirac 
operators  on non--compact
 almost complex, complete orbifolds which are sufficiently
 regular at infinity, admit a unique extension.
 Additonally, we prove a  generalized orbifold Stokes'/Divergence 
 theorem.
 \endabstract \endtopmatter

					\noindent {\bf 0.  Introduction.}
					 
					  Orbifolds, generalized manifolds that are locally
					  the quotient of an euclidean  space modulo a
					  finite group of isometries, were first introduced first by
					  Satake. In the late seventies,  Kawasaki  proved an
					  orbifold signature formula, together with more
					  general index theorems, see [Kw1], [Kw2], [Kw3]. In
					  [Fa1] we  proved a $K$--theoretical index theorem
					  for orbifolds with operator algebraic means,
					  and  in [Fa2], [Fa3] we studied compact
					   orbifold spectral theory and defined orbifold eta invariants. Other orbifold index formulas were proved in [Du], [V].  In [Ch] Chiang studied compact
					   orbifold heat kernels and harmonic maps, while in [Stan], Stanhope established some interesting
					   geometrical applications of orbifold spectral
					   theory.   

					   Here we will continue the orbifold spectral
					   analysis started in [Fa2] and [Fa3]. In
					   particular we show that on a non--compact complete almost complex $Spin^c$ orbifold which is sufficiently regular at
					   infinity (see Definition 2.1), generalized Dirac operators are closed. This extends to
					   orbifolds theorems of Gaffney [Gn1], Yau [Y],
					   and Wolf [W], whose ideas are used in our proofs,  together with more orbifold-specific techniques.
					   In particular, our first main result, Theorem
					   3.1, asserts that

					    \proclaim{Theorem 3.1 } Let $X$ be an even--dimensional non--compact complete Hermitian 
$Spin^c$ almost complex orbifold which is
sufficiently regular at infinity. Assume that a Hermitian connection is chosen on the dual of its canonical line bundle $K^*$. Let $E$ be a proper Hermitian orbibundle (with connection $\nabla^E$) over
  $X$, and let $D_E$ be the generalized Dirac operator with coefficients in $E$. 
								   Let ${\Cal D}( D_E^{MIN})$ be the
								   domain of the min extension of $D_E$,
								   and  $ {\Cal D}( D_E^{MAX})  $ be the
								   domain of the max extension of $D_E$,
								   see the end of Section 2 for details.
								   Then
								   $$
								   {\Cal D}( D_E^{MIN})= {\Cal D}(
								   D_E^{MAX}).  
								   $$
								   \endproclaim

					   We also prove the following
					   Stokes'/Divergence theorem, Theorem 5.1, which
 generalizes to  orbifolds results of Gaffney
					   [Gn2], Karp [K],  and Yau [Y].
						
\proclaim{Theorem 5.1} Let $X$ be an even--dimensional  non--compact complete $Spin^c$ almost complex orbifold which is sufficiently regular at infinity.
Assume that a connection is chosen on the dual of its canonical line bundle. Let $V$ be a vector field on $X$ such that 
$$
\lim_{k\to +\infty} \inf
\frac1k \int_{B_{2k}-B_k}
\Vert V\Vert\, dv =0,
 $$
where $\Vert V\Vert$
denotes the length of $V$, and $B_k= \{ y\in X|
\rho(y) =d(y,y_0)\leq k\}$ for a fixed $y_0 \in X-\Sigma(X)$, where $\Sigma(X)$ is the singular locus of $X$. 
 Then if either $(div\,
V)^+$ or $(div\, V)^-$ is
integrable on $X$, we have
$$
\int_X div\,(V) \, dv =0.
$$ 
\endproclaim

						In a sequel to this paper [Fa4], we 
						use the results we proved here to establish an orbifold Gromov-Lawson relative index theorem,
						c.f. [GL] for the manifold case.  
						More in detail, the contents of this paper are as follows.  In Section 1
						we recall the definition of orbifold, orbibundles, and introduce orbifold Dirac operators.
						In Section 2, we study Dirac operators on non--compact
						orbifolds 
						from a local viewpoint. In Section 3, we state
						and prove our first main result, Theorem 3.1.
						In Section 4 we prove that, if $D$ is a Dirac
						operator, 
						and $D^2\sigma=0$, also $D\sigma=0$. In Section	5 we finally prove our Stokes'/Divergence
theorem, Theorem 5.1. Vanishing results are considered in Section 6.

						In the sequel,  all orbifolds and manifolds are even--dimensional, smooth, Hermitian,						$Spin^c$, connected, and almost complex 
unless otherwise specified. All vector and orbibundles are assumed to be smooth and proper. We also assume that all of our orbifolds/manifolds are endowed with a fixed Hermitian connection on the dual of their canonical line bundle $K^*$. This latter hypotheses  
allows us to define a \lq canonical\rq $\,Spin^c\,$ 
Dirac operator and, given a Hermitian orbibundle $E$ with a chosen connection $\nabla^E$, the \lq canonical\rq $\, Spin^c \,$ 
Dirac operator with coefficients in $E$. Both of these operators depend, in the $Spin^c$ case, on the choice of the selected connections, see [Du; Chapter 14], and [LM; Appendix D]. For the $Spin$ or {\it complex} 
case, the choice of the connection on $K^*$ is canonical.

I would like to thank the sabbatical program of the University of Colorado/Boulder, and the Mathematics Department of  the University of Florence,
Italy, for their warm hospitality during the period this paper was written. We also thanks the referee for useful suggestions. 

\vskip 1em
\noindent {\bf 1. Orbifolds, Orbibundles and Dirac Operators.} \vskip 1em
							
							 In this section we will review some
							 definitions and results that we will use  throughtout this paper. For generalities on orbifolds and operators on orbifolds, see [Kw1], [Kw2], [Kw3], [Ch], [Du], [V]. 

							 An orbifold is a Hausdorff second countable topological space $X$ together with an
							 atlas of charts 
							 $\Cal U  =\{ {(\tilde U_i}, G_i )| i \in { I\} }$, with ${\tilde U}_i/G_i = U_i$ open,  and
							with projection $\pi_i: {\tilde U}_i \to U_i$, $i\in I$, 
							 satisfying the following properties 
							 \roster
							 \item  If  two charts $U_1$ and $U_2$
							 associated to pairs $({\tilde U}_1, G_1)$, $({\tilde U}_2, G_2)$ of $\Cal U$,  are such that $U_1\subseteq U_2$, then there exists a smooth
							 open embedding $\lambda$:  ${\tilde U}_1$
							 $\to {\tilde U}_2$ and a homomorphism
							 $\mu$: $G_1$ $\to G_{2}$ such that $\pi_1=
							 \pi_{2}\circ \lambda$ and $\lambda \circ
							 \gamma =\mu(\gamma)\circ \lambda$, $\forall
							 \gamma \in G_1$.
							 \item The collection of the open charts
							 $U_i$, $i\in I$, belonging to the atlas
							 $\Cal U$ forms a basis for the topology on $X$. 
							 \endroster
							 We will call an orbifold atlas as above a
							 standard orbifold atlas.

							 For any $x$ point of $X$, the isotropy
							 $G_x$ of $x$ is well defined, up to
							 conjugacy, by using any local coordinate
							 chart. The set of all points $x\in X$ with
							 non--trivial stabilizer, $\Sigma(X)$, is
							 called the singular locus of $X$, see e.g.
							 [Ch].  Note that $X-\Sigma(X)$ is a smooth
							 manifold.
								
								If we now endow $X$ with a countable
								locally finite orbifold atlas $\Cal F$, $\Cal F  =\{(
								{\tilde U_i}, G_i )| i \in {\Bbb N\} }$,
							then 	by standard theory
								there exists a smooth partition of unity
								$\eta=\{ \eta_i \}_{i\in \Bbb N}$
								subordinated to $\Cal F$, [Ch]. This in
								particular means that, for any $i \in
								\Bbb N$, $\eta_i$ is a smooth function
								on $U_i$ (i.e., its lift to any chart of
								a standard orbifold atlas is smooth),
								the support of $\eta_i$ is included in
								an open subset  $U_i'$ of $U_i$, and $\cup U_i'= X$. We will call any $\eta$ as above an $\Cal F$-partition of unity. 
								Let $E$ be an orbibundle
								over the orbifold $X$. (For the precise
								definition see [Kw1], [Kw2], [Kw3],
								[Ch].) In particular $E$ is an orbifold
								in its own right; on an orbifold chart
								$U_1$ associated to a pair $({\tilde
								U}_1, G_1)$ of a standard orbifold atlas
								$\Cal U  =\{ ({\tilde U_i}, G_i )| i \in {
								I\} }$ of $X$, $E$ lifts to a
								$G_1$-equivariant bundle. Standard
								orbifold atlases on $X$ can be used
								to provide standard orbifold atlases on
								$E$. 

								If $E$ is an  orbibundle
								over the orbifold $X$, a section $s:X\to
								E$ is called a smooth orbifold section
								if for each chart $U_i$ associated to a
								pair $({\tilde U}_i, G_i)$ of a standard
								orbifold atlas $\Cal U  =\{ ({\tilde
								U_i}, G_i) | i \in { I\} }$ of $X$, we
								have that $s|_{U_i}: U_i \to E|_{U_i}$
								is covered by a smooth $G_i$--invariant
								section ${\tilde s}|_{{\tilde U}_i}:
								{\tilde U}_i \to {\tilde E|_{{\tilde U}_i}}$.
								Given a Hermitian orbibundle $E$ over $X$,
								we will denote by 
								${\Cal C}^\infty(X,E)$ the space of all
								smooth sections of $E$, and by
								${\Cal C}^\infty_c(X,E)$  the space of
								all smooth sections with compact
								support.
								Classical orbibundles over $X$
								are the tangent bundle $TX$, and the
								cotangent bundle $T^*X$ of $X$. We can
								form orbibundles tensor products by
								taking the tensor products of their local
								expressions in the charts of a standard
								orbifold atlas.  

								Define an inner product between sections 
								of ${\Cal C}^\infty (X,E)$  (or
			${\Cal C}^\infty_c(X,E)$) by a the
	following formula (c.f., [Ch; 2.2a])
								$$
								(\sigma_1, \sigma_2)=
								\sum_{i=1}^{+\infty} \frac1{|G_i|}
								\int_{\tilde U_i}
								 {\tilde \eta_i} (\tilde x_i ) <{\tilde
								 \sigma_1}(\tilde x_i ), {\tilde
								 \sigma_2}(\tilde x_i )> dv(\tilde x_i
								 ),
								 $$ 
								 where $\eta=\{ \eta_i \}_{i\in \Bbb
								 N}$, is a $\Cal F$-partition of unity
								 subordinated to the locally finite
								 orbifold cover $\Cal F  =\{ ({\tilde
								 U_i}, G_i) | i \in {\Bbb N} \}$, and $<,>$ is a $G_i$--invariant product on $\tilde E$. (Note that, with slight abuse of
								 notation, we used $\,\tilde {}\,$ to denote
								 lift to ${\tilde U_i}$.)

								 We will now review the construction of the 
 Dirac operators with coefficients in an Hermitian 
								 orbibundle $E$ endowed with a connection $\nabla^E$, over an 
orbifold $X$ satisfying our hypotheses, [Du; Sections 5 and 12], [Kw2], [BGV], [LM; Appendix D]. First of all, $X$  admits a
								 $Spin^c$-principal tangent orbibundle,
								 $Spin^c(TX)$, with, in our hypotheses,  canonical $Spin^c$
								 orbifold connection $\nabla^c$.  Let $\Delta^{\pm,
								 c}$ be the half  $Spin^c$
								 representations (recall that the $X$ is even dimensional). Then we have two
								 orbibundles
								 $$
								 \Delta^{\pm, c} (TX) = Spin^c(TX)
								 \times_{ Spin^c} \Delta^{\pm, c}, $$
								 with induced connections $\nabla^{\pm,
								 c}$, from $\nabla^c$; $\nabla^{\pm,c} : {\Cal C}^\infty_c (X, \Delta^{\mp, c} (TX))
\to {\Cal C}^\infty_c (X, T^*X \otimes \Delta^{\mp, c} (TX))$.								 The Clifford module structure on
								 $\Delta^{\pm, c}$ defines Clifford
								 multiplications
								 $$
								 m_{\pm}: TX \otimes_{\Bbb R}
								 \Delta^{\pm, c} (TX) \to \Delta^{\mp,
								 c} (TX)
								 $$ 
								 On $E$ we have the  
								 connection $\nabla^E$. Then the
								 generalized $\pm $ Dirac operator with
								 coefficients in $E$, $d_E^{\pm, c}$, 
								 $$
								 d_E^{\pm, c}: {\Cal C}^\infty_c (X,
								 \Delta^{\pm, c} (TX)\otimes_{\Bbb C} E)
								 \to
								 {\Cal C}^\infty_c (X, \Delta^{\mp, c}
								 (TX)\otimes_{\Bbb C} E)
								 $$
								 is defined by
								 $$
								 d_E^{\pm, c}= M \circ \left( \nabla^{\pm, c}\otimes Id + Id \otimes \nabla^{E}\right),
								 $$
								 where $M$ denotes the map induced by Clifford multiplication and TX has been identified with $T^*X$ via
								 the orbifold  metric. 
We will also use the notation  $\Cal S$ for
								 the orbifold $Spin^c$ bundle
								 $(\Delta^{+, c}\oplus \Delta^{-, c})	 (TX)$, and ${\Cal S}\otimes E$ or $\Cal E$  for $(\Delta^{+,
								 c}\oplus \Delta^{-, c})
								 (TX)\otimes_{\Bbb C} E$ throughout this paper. 
 We will define $D_E$, the generalized
								 Dirac operator on $X$ with coefficient
								 in $E$, to be $(d_E^{+, c}+ d_E^{-, c})$.

								 \vskip 1em
								 \noindent {\bf 2.  Dirac Operators on
								 Complete Orbifolds.}

								 \vskip 1em

								 On an orbifold $X$ (not necessarily
								 compact), the generalized Dirac
								 operator with coefficients in the orbibundle $E$ (with connection $\nabla^E$), $D_E$, as
								 defined in Section 1, is given by 
								 $$
								 D_E: {\Cal C}^\infty_c (X, (\Delta^{+,
								 c}\oplus \Delta^{-, c})
								 (TX)\otimes_{\Bbb C} E) \to {\Cal
								 C}^\infty_c (X, (\Delta^{-, c}\oplus
								 \Delta^{+, c}) (TX)\otimes_{\Bbb C} E)
								 $$
								 $$
								 D_E= M \circ \left( (\nabla^{+, c} +
								 \nabla^{-, c})\otimes Id + Id \otimes
								 \nabla^{E} \right).
								 $$

								 On orbifold charts, the Dirac operator
								 $D_E$ with coefficients in the 
								 Hermitian orbibundle $E$ (with connection $\nabla^E$), has the following local
								 expression ${\tilde D}_E$. Let $\Cal U
								 =\{({\tilde U_i}, G_i) | i \in { I\} }$,
								 with ${\tilde U}_i/G_i = U_i$ be a
								 standard orbifold atlas. On a local
								 chart ${\tilde U}_i$, $i\in I$ fixed, we have
								 $$
								 \Delta^{\pm, c} (T{\tilde U_i}) =
								 Spin^c(T{\tilde U_i}) \times_{ Spin^c}
								 \Delta^{\pm, c}, $$
								 with induced $G_i$--invariant
								 connections $\nabla^{\pm, c}$, from
								 $\nabla^c$.
								 The Clifford module structure on
								 $\Delta^{\pm, c}$ defines Clifford
								 multiplications
								 $$
								 m_{\pm}: T{\tilde U_i} \otimes_{\Bbb R}
								 \Delta^{\pm, c} (T{\tilde U_i}) \to
								 \Delta^{\mp, c} (T{\tilde U_i}).
								 $$ 
								 On $\tilde E$, the lift of $E$, we have the 
								 $G_i$--invariant connection
								 $\nabla^{\tilde E}$. 
								 Then the generalized $\pm $
								 Dirac operators with coefficients in
								 $E$, ${\tilde d}_E^{\pm, c}$, 
								 $$
								 {\tilde d}_E^{\pm, c}: {\Cal
								 C}^\infty_c ({\tilde U_i}, \Delta^{\pm,
								 c} (T{\tilde U_i})\otimes_{\Bbb C}
								 {\tilde E}) \to
								 {\Cal C}^\infty_c ({\tilde U_i},
								 \Delta^{\mp, c} (T{\tilde
								 U_i})\otimes_{\Bbb C} E)
								 $$
								 is given by
								 $$
								 {\tilde d}_E^{\pm, c}= M \circ \left(
								 \nabla^{\pm, c}\otimes Id + Id
								 \otimes \nabla^
								 {\tilde E} \right),
								 $$
								 where $M$ is induced by Clifford multiplication and $T{\tilde U_i}$ has been identified
								 with $T^*{\tilde U_i}$ via the
								 $G_i$--invariant metric. 
								 Also, ${\tilde D}_E$, the generalized
								 Dirac operator on $X$ with coefficient
								 in $E$, is given by ${\tilde d}_E^{+,
								 c}+ {\tilde d}_E^{-, c}$ on ${\tilde U}_i$. 

								 If $e_1, \dots, e_n$ is an orthonormal
								 local basis for the space $T{\tilde
								 U_i}$ at a point $\tilde x$, then ${\tilde
								 D}_E$
								 has local expression
								 $$
								 {\tilde D}^{ E} =
								 \sum_{k=1}^{n} e_k {\tilde
								 \nabla}^{E}_{e_k},
								 $$

								 where 
								 $$
								 {\tilde \nabla}^{E}= ({\nabla}^{+, c}+
								 {\nabla}^{-, c} ) \otimes 1 + 1 \otimes
								 {\nabla}^{\tilde E}.
								 $$

								 Now, in analogy with the manifold case,
								 see [GL], [W], [Gn1], [LM], 
[Y], we
								 will show that $D_E$ is symmetric,
								 whenever $X$ is a sufficiently regular at
								 infinity.

\proclaim{Definition 2.1 } Let $X$ be a
								 non--compact complete 
orbifold. Then we say that
								 $X$ is sufficiently regular at infinity
								 if, for any neighborhood  $\Omega
								 \subseteq X$ of infinity, there exists
a compact domain $K_\Omega$  with $\Omega \cup K_\Omega =X$ and with  boundary  strictly included in $\Omega$,
								 on which the Divergence and Stokes'
								 Theorems hold.  
								 \endproclaim

								 For a compact orbifold
								 without boundary, the Divergence
								 Theorem holds,
								 [Ch]. See also [C] for other
								 results. Sufficient regularity also holds in
								 the case of a product end, by an
								 adaptation of Chiang's method, 
								 [Ch], and in the case of finite volume
								 hyperbolic orbifolds because of the structure 
								 of the cusps cross sections, [LoR]. Also, 
								 geometrically finite orbifolds with pinched negative sectional curvature
								 satisfy this hypothesis, [AX]. In general, ours seems to be a
								 very reasonable assumption to make,
								 which will be certainly satisfied in
								 many cases of interest, see above examples. 
								 For Sobolev inequalities of Gallot type involving domains, see [N].
								 \proclaim{Theorem 2.2 } Let $X$ be a
							non--compact complete orbifold which is
								 sufficiently regular at infinity, and
								 let $E$ be a Hermitian orbibundle (with connection $\nabla^E$) over $X$.
								 Let $D_E$ be the generalized
								 Dirac operator with coefficients in
								 $E$, as defined above.
Then $D_E$ is symmetric, i.e.,    
								 $$
								 (D_E \sigma_1, \sigma_2)= (\sigma_1,
								 D_E\sigma_2), \quad \forall \sigma_1,
								 \sigma_2
								 \in {\Cal C}^\infty_c (X,
								 S\otimes_{\Bbb C} E),
								 $$
								 where $(,)$ denotes the inner product
								 defined earlier.
								 \endproclaim

								 \noindent {\it Proof.} Let $\Cal E=
								 S\otimes E$, $D=D_E$.   
								 Let $\eta=\{ \eta_i \}_{i\in \Bbb N}$
								 be a $\Cal F$-partition of unity
								 subordinated to the locally finite
								 orbifold cover $\Cal F  =\{ (\tilde
								 U_i, G_i )| i \in {\Bbb N} \}$. (Note
								 that, as before, we are using $\, \tilde
								 {}\,$ to denote  lift to ${\tilde U_i}$.)
								 Then
								 $$
								 (\sigma_1, \sigma_2)=
								 \sum_{i=1}^{+\infty} \frac1{|G_i|}
								 \int_{\tilde U_i}
								  {\tilde \eta_i} (\tilde x_i ) <{\tilde
								  \sigma_1}(\tilde x_i ), {\tilde
								  \sigma_2}(\tilde x_i )> dv(\tilde x_i
								  ), \quad \forall \sigma_1, \sigma_2
								  \in
								  {\Cal C}^\infty_c ({U}_i, {\Cal
								  E}).
								  $$ 
								  Since $T({\tilde U}_i )$ is
								  parallelizable, we can choose a local
								  orthonormal basis $e_1, \dots, e_n$
								  for the space $T{\tilde U_i}$ at any
								  point ${\tilde x}_i$; thus, if we set
								  $\nabla = \tilde {\nabla}^E$, 
								  $$
								  <{\tilde D}, {\tilde \sigma}_1
								  ({\tilde x}_i), {\tilde
								  \sigma}_2({\tilde x}_i)>=
								  \sum_{k=1}^n <e_k \dot \nabla_{e_k}
								  {\tilde \sigma}_1 ({\tilde x}_i),
								  {\tilde \sigma}_2 ({\tilde x}_i)>    
								  $$
								  $$
								  -\sum_{k=1}^n <\nabla_{e_k} {\tilde
								  \sigma}_1 ({\tilde x}_i), e_k \dot
								  {\tilde \sigma}_2 ({\tilde x}_i)>
								  $$
								  $$
								  -\sum_{k=1}^n \left\{ \nabla_{e_k}
								  <{\tilde \sigma}_1 ({\tilde x}_i), e_k
								  \dot {\tilde \sigma}_2 ({\tilde x}_i)>
								  -  <{\tilde \sigma}_1 ({\tilde x}_i),
								  (\nabla_{e_k} )\dot {\tilde \sigma_2
								  ({\tilde x}_i) + e_k \dot \nabla_{e_k}
								  \sigma}_2 ({\tilde x}_i)> \right\},
								  $$
								  where  $<,>$ is a $G_i$--invariant
								  inner product on $\tilde {\Cal S}$.

								  If we define the $G_i$--invariant vector field $V_i$ on ${\tilde U}_i$ by
								  $$
								  <V_i, W > =-  <{\tilde \sigma}_1, W
								  \circ {\tilde \sigma}_2  >, \, \hbox{ for any vector field } W,
								  $$
								  we have that the above expression can  be rewritten as
								  $$
<{\tilde D}, {\tilde \sigma}_1 ({\tilde x}_i), {\tilde \sigma}_2({\tilde x}_i)>
 =div (V_i({\tilde x}_i)) + <{\tilde \sigma}_1 ({\tilde x}_i), <{\tilde D}{\tilde
								  \sigma}_2 ({\tilde x}_i)> \tag 2.1
								  $$
								  Now integrate (2.1) (multiplied by
								  ${\tilde \eta}_i$ and divided by
								  $|G_i|$)
								  over ${\tilde U}_i$.  Then by using the
								  Divergence Theorem, we are done.   $\qed$

								  \proclaim{Remark 2.3 } Theorem 2.2 is
								  also valid when only one of the two
								  sections $\sigma_1$, $\sigma_2$ has
								  compact support.
								   \endproclaim

								   Now complete the space ${\Cal
								   C}^\infty_c (X, {\Cal E})$, ${\Cal
								   E}= S\otimes_{\Bbb C} E$, $\Cal S$
								   $Spin^c$ bundle on $X$, $E$ Hermitian 
								   orbibundle over $X$, with respect to
								   the norm 
								   $$
								   \Vert \sigma \Vert_X= \sqrt{ <\sigma,
								   \sigma>}=
								   \left(\sum_{i=1}^{+\infty}
								   \frac1{|G_i|} \int_{\tilde U_i}
								   {\tilde \eta_i} (\tilde x_i )
								   <{\tilde \sigma_1}(\tilde x_i ),
								   {\tilde \sigma_2}(\tilde x_i )>
								   dv(\tilde x_i )
\right)^{\frac12}. 
								   $$
								  We thus obtain the
								   ${\Cal L}^2$--space ${\Cal L}^2(X,
								   {\Cal E})$.  The Dirac operator
								   $$
								   D_E: {\Cal C}^\infty_c (X, {\Cal E})
								   \to {\Cal C}^\infty_c (X, {\Cal E}) 
								   $$
								   has two natural extensions, min and max listed below,  as an
								   unbounded operator
								   $$
								   D_E: {\Cal L}^2 (X, {\Cal E}) \to
								   {\Cal L}^2 (X, {\Cal E}). 
								   $$

								   1. Minimal Extension $D_E^{MIN}$.
								   The minimal extension of $D_E$,
								   $D_E^{MIN}$, is obtained by taking
								   the graph closure of the graph of
								   $D_E$, i.e.,
								   $$
								   D_E^{MIN}: {\Cal D}( D_E^{MIN}) \to
								   {\Cal L}^2 (X, {\Cal E}), 
								   $$
								   where ${\Cal D}( D_E^{MIN})$, the domain of $D_E^{MIN}$, is defined to
								   be   
								   the set of $\sigma\in {\Cal L}^2 (X,
								   {\Cal E})$ for which there exists a
								   sequence $\sigma_k\in {\Cal
								   C}^\infty_c (X, {\Cal E})$ such that
								   $\sigma_k \to \sigma$ and
								   $D_E\sigma_k \to \tau$ in  ${\Cal
								   L}^2 (X, {\Cal E})$,
								   for some $\sigma, \tau \in {\Cal L}^2
								   (X, {\Cal E})$. Set
								   $D_E^{MIN}(\sigma)=\tau$.

								   2. Maximal Extension $D_E^{MAX}$.
								   The maximal extension of $D_E$,
								   $D_E^{MAX}$, is obtained by taking
								   its domain to be the set of all
								   $\sigma \in {\Cal L}^2 (X, {\Cal E})$
								   such that the distributional image of
								   $D_E(\sigma)$ is still in ${\Cal L}^2
								   (X, {\Cal E})$. More precisely,
								   $$
								   D_E^{MAX}: {\Cal D}( D_E^{MAX}) \to
								   {\Cal L}^2 (X, {\Cal E}), 
								   $$
								   where ${\Cal D}( D_E^{MAX})$, the
								   domain of $D_E^{MAX}$, is defined to
								   be   
								   the set of $\sigma\in {\Cal L}^2 (X,
								   {\Cal E})$ such that the linear
								   functional $L(\sigma_2)=(\sigma,
								   D_E(\sigma_2))$ on  $ {\Cal
								   C}^\infty_c (X, {\Cal E})$ is bounded
								   in the ${\Cal L}^2 (X, {\Cal E})$
								   norm. Note that the boundedness of
								   $L$ implies that there exists an
								   element $\tau \in {\Cal L}^2 (X, {\Cal
								   E})$, such that
								   $$
								   (\tau, \sigma_2)= (\sigma,
								   D_E\sigma_2), \quad \forall 
								   \sigma_2 \in {\Cal C}^\infty_c (X,
								   {\Cal E}).
								   $$  
								   Define the above $\tau$ to be
								   $D_E^{MAX}(\sigma)$.

								   \proclaim{Remark 2.4 } Since $( \, ,\, 
								  )$ is continuous in the ${\Cal
								   L}^2 (X, {\Cal E})$ norm,  we have
								   that 
								   $$
								   {\Cal D}( D_E^{MIN}) \subseteq {\Cal
								   D}( D_E^{MAX}) . 
								   $$
								   \endproclaim

								   \vskip 1em
								   \noindent {\bf 3.  Generalized Dirac
								   Operators on Non--Compact
Complete Orbifolds are
								   Closed.}
								   \vskip 1em

								   In this Section we will prove that
								   generalized Dirac operators on
								   complete orbifolds which are
								   sufficiently regular at infinity, are
								   closed operators. This theorem
   generalizes to orbifolds [GL; Theorem 1.17], and [W; Theorem 5.1].

\proclaim{Theorem 3.1 } Let $X$ be an even--dimensional non--compact complete Hermitian 
$Spin^c$ almost complex orbifold which is
sufficiently regular at infinity. Assume that a Hermitian connection is chosen on the dual of its canonical line bundle $K^*$. Let $E$ be a proper Hermitian orbibundle (with connection $\nabla^E$) over
  $X$, and let $D_E$ be the generalized Dirac operator with coefficients in $E$. 
								   Let ${\Cal D}( D_E^{MIN})$ be the
								   domain of the min extension of $D_E$,
								   and  $ {\Cal D}( D_E^{MAX})  $ be the
								   domain of the max extension of $D_E$,
								   see the end of Section 2 for details.
								   Then
								   $$
								   {\Cal D}( D_E^{MIN})= {\Cal D}(
								   D_E^{MAX}).  
								   $$
								   \endproclaim

								   Our proof of Theorem 3.1, which will
								   occupy the remaining of this section, 								   will be an adaptation of [W; Proof
								   of Theorem 5.1]. In particular,
								   suitable modifications to Wolf's
								   proof for manifolds will be mostly
needed to deal  
with orbifold distance functions.

\noindent {\it Proof.}  Firstly,
								   recall that we denoted by $\Sigma(X)$   the singular locus of $X$.
								   Then $X-\Sigma(X)$ is a convex						   manifold. In particular any two
							   points of $X-\Sigma(X)$ can be
								   connected by a geodesic arc lying
								   entirely in $X-\Sigma(X)$, see
								   [Stan; Section 4].     Because of Remark 2.4, to prove
								   Theorem 3.1 it is enough to show that 
								   $$
								   {\Cal D}( D_E^{MAX})\subseteq {\Cal
								   D}( D_E^{MIN}). \tag 3.1 
								   $$
								   Note that ${\Cal D}( D_E^{MAX})$
								   carries the norm
								   $$
								   N(\sigma)=\left\{  \Vert \sigma
								   \Vert^2_X +\Vert D_E(\sigma) \Vert^2_X
								   \right\}^{\frac12}, \quad \forall
								   \sigma \in {\Cal D}( D_E^{MAX}),
								   $$
								   where $\Vert\, \Vert_Y$ denotes the
								   ${\Cal L}^2 (Y, {\Cal E})$ norm, for
								   $Y\subseteq X$, c.f. Section 2. 
								   But (3.1) is equivalent to 
								   $${\Cal C}^\infty_c (X, {\Cal E})
								   \hbox { is dense in }{\Cal D}(
								   D_E^{MAX})  \hbox{  in the norm N}.
								   \tag 3.2 $$ 
								   Thus Theorem 3.1 will clearly follow
								   once we have proven the Lemmas 3.2
								   and 3.3 below. $\qed$.

								   \proclaim{Lemma 3.2 } Let $X$ be a non--compact
								   complete orbifold which is
								   sufficiently regular at infinity, and
								   let $E$ be a Hermitian orbibundle (with connection $\nabla^E$)  over
								   $X$. Let
								   $D_E: {\Cal C}^\infty_c (X,
								{\Cal E}) \to {\Cal
C}^\infty_c (X, {\Cal E}),$
								   be the generalized Dirac operator on
								   $X$ with coefficients in $E$ as in
								   Theorem 3.1. Then
								   $$
								   {\Cal C}^\infty_c (X, {\Cal E})
								   \hbox{ is dense in }{\Cal D}_c(
								   D_E^{MAX})  \hbox{  in the norm N},
								   $$
where we set ${\Cal D}_c( D_E^{MAX}) $ to be the subset of the elements of
								   ${\Cal D}( D_E^{MAX})$ with compact
								   support, and where 
								   $$
								   N(\sigma)=\left\{  \Vert \sigma
								   \Vert^2_X +\Vert D_E(\sigma) \Vert^2_X
								   \right\}^{\frac12}, \quad \forall
								   \sigma \in {\Cal D}( D_E^{MAX}).
								   $$
								     \endproclaim

								   \proclaim{Lemma 3.3 } Let $X$ be a non--compact complete orbifold which  is
sufficiently regular at infinity, and
let $E$ be a Hermitian orbibundle (with connection $\nabla^E$)  over
$X$. Let
$D_E: {\Cal C}^\infty_c (X,
{\Cal E}) \to {\Cal
C}^\infty_c (X, {\Cal E})
$ be the generalized Dirac operator on
								   $X$ with coefficients in $E$ as in Theorem 3.1 and Lemma 3.2. Then
$$
								   {\Cal D}_c( D_E^{MAX})  \hbox{ is
								   dense in }{\Cal D}( D_E^{MAX})
								   \hbox{  in the norm N}. $$ 
								   \endproclaim

								   \noindent {\it Proof of Lemma
3.2.}
								   Set ${\Cal D}_c = {\Cal D}_c(
								   D_E^{MAX}) $, and 
let   $\sigma \in {\Cal
								   D}_c$. Choose a locally finite
								   orbifold atlas $\Cal F$,  $\Cal F
								   =\{ ( {\tilde U_i}, G_i ) | i \in {\Bbb
								   N\} }$, with associated smooth
								   partition of unity $\eta=\{ \eta_i
								   \}_{i\in \Bbb N}$. Suppose that 
								   $supp(\eta_i)\cap supp (
								   \sigma)\not=\emptyset$ only for $i=1,
								   \dots, \ell$. Then
								   $\sigma= \sigma_1+\dots,
								   \sigma_\ell$, with $\sigma_i= \eta_i
								   \sigma$ having support in $U_i$,
								   $i=1,\dots, \ell$. We can lift
								   $\sigma_i$ to a $G_i$--invariant
								   section $\tilde {\sigma}_i$. By
								   trivializing the bundle $\tilde E$
								   over ${\tilde U}_i $, we can assume
								   that we are dealing with functions.
								   Convolutions with an approximated
								   identity and averaging, give a
								   $G_i$--invariant sequence $\{{\tilde
								   u}_{i,k}\}_{k\in N}$ in $ {\Cal
								   C}^\infty_c({\tilde U}_i, {\tilde
								   {\Cal E}})$ whose image $\{{u}_{i,k}\}_{k\in N}$
								   in   $ {\Cal C}^\infty_c({U}_i,
								   {{\Cal E}})$ satisfies $N(\sigma_i -
								   {u}_{i,k}) < \frac 1k$, $i=1,\dots,
								   \ell$.  (The $\Cal L^2$--norm is
								   computed by dividing by $|G_i|$ and
								   integrating on ${\tilde U}_i$.) Now,
								   if we set, $u_k= {u}_{1,k}+\dots +
								   {u}_{\ell,k}$, we have
								   $N(\sigma- u_k)< \frac{\ell}k$, so
								   ${u}_{k} \to \sigma$ in the norm N.
									 $\qed$

									 \noindent {\it Proof of Lemma 3.3.}
									 Set ${\Cal D}_c = {\Cal D}_c(
									 D_E^{MAX}) $, ${\Cal D} = {\Cal
									 D}( D_E^{MAX}) $, and $D=D_E$. Let $y_0\in
									 X-\Sigma(X)$ be fixed, where
									 $\Sigma(X)$ is the singular locus
									 of $X$. Let $y\in X$ and denote by
									 $\rho(y)$ the orbifold distance
 bewteen $y_0$ and $y$. We will only
									 be interested in the beahviour of
									 $\rho$ at points of the convex
									 manifold $X-\Sigma(X)$. (For more
									 details on $X-\Sigma(X)$, see
									 [Stan] and [B].) Note also that
									 $\Sigma(X)$  has measure zero in
									 $X$, [Ch], and so $\rho$ is a
									 function which is differentiable on
									 $X$ a. e.  Therefore
									 we can assume that we have
									 $$
									 \Vert \nabla\rho \Vert \leq 1 \quad \hbox
		 { almost everywhere on X,}
									 $$   
where the above norm is the sup norm,   
									 c.f. [W; pg. 623]; we need to
									 additionally remove the measure
									 zero set $\Sigma(X)$. We can now
									 proceed as in  [W; Proof of
									 (5.5)]. For completeness, we go
									 through all the details of the
									 proof below.
									 If $r>0$, let 
									 $$
									 B_r= \left\{  y\in X\, | \,\rho(y)<r
									 \right\}.
									 $$
									 Since $X$ is complete the closure
									 of $B_r$, $\overline{ B}_r$, is
									 compact.
									 Choose a ${\Cal C}^\infty$ function
									 $a: \Bbb R \to [0,1]$ such that
									 $a(-\infty, 1]=1$, $a[2,
									 +\infty)=0$, and denote by $M$ the
									 max of $a'$ on $\Bbb R$. If $r>0$
									 as before, define
									 $$
									 b_r :X\to [0,1], \quad \quad\hbox { by }\quad
									 b_r(y) = a (\frac{\rho(y)}r).
									 $$ 
									 Then
									 $$
									 b_r =1\hbox { on } B_r, \quad
									 supp(b_r)\subseteq
									 \overline{B}_{2r}.
									 $$
									 We have that $b_r$ is
 differentiable almost everywhere,
									 and, at points of
									 differentiability, the following inequality holds
									 $$
									 \Vert \nabla (b_r) \Vert^2 = \Vert
									 \frac1r a'(\frac{\rho}r)\Vert^2
									 \leq \frac{M^2}{r^2}.
									 $$
									 Fix $\sigma \in \Cal D$, and  write
									 $\sigma_s = b_s \sigma$, for $s\in
									 \Bbb N$. Now 
									 $\sigma_s\in {\Cal D}_c $, since
									 the support of $b_s$ is contained
									 in $\overline{B}_{2s}$ compact.  
									 Choose a locally finite orbifold
									 atlas $\Cal F$,  $\Cal F  =\{
									 ({\tilde U_i}, G_i ) | i \in {\Bbb N\}
									 }$, with associated smooth
									 partition of unity $\eta=\{ \eta_i
									 \}_{i\in \Bbb N}$. Suppose that
									 $supp(\eta_i)\cap
									 \overline{B}_{2k}\not=\emptyset$ only
									 for $i=1, \dots, \ell$. Then on a
									 local chart ${\tilde U_i}$, $ =1,
									 \dots, \ell$, we have, as in the
									 proof of Theorem 2.2 (as usual denote
									 by $\, \tilde {}\,  $ the lift to
									 ${\tilde U}_i$),  
									 $$
									 {\tilde D} ( {\tilde \sigma}_s )=
									 {\tilde D} ( {\tilde b}_s {\tilde
									 \sigma} )=
									 \sum_{j=1}^n e_j\left(
									 \nabla_{e_j} ({\tilde b}_s {\tilde
									 \sigma}) \right)
									 $$
									 $$
									 =\sum_{j=1}^n e_j\left( e_j
									 ({\tilde b}_s) {\tilde \sigma} +
									 ({\tilde b}_s \nabla_{e_j} {\tilde
									 \sigma}) \right)
									 $$
									 $$
									 = \nabla  ({\tilde b}_s) {\tilde
									 \sigma}+ ({\tilde b}_s) {\tilde
									 D}(\tilde \sigma)
									 $$
									 almost everywhere on ${\tilde
									 U}_i$. Since $b_s=1$ on $B_s$, we
									 have
									 $$
									 \Vert D(\sigma -\sigma_s)
									 \Vert^2_X = \Vert (1-b_s)
									 D(\sigma) +\nabla (b_s) \sigma
									 \Vert^2_X  									 $$
									 $$
									 \leq \Vert D(\sigma)
									 \Vert^2_{X-B_s}+ \frac{M^2}{s^2}
									 \Vert \sigma \Vert^2_X,
									 $$
									 where $\Vert \, \Vert_Y$ denotes
									 the ${\Cal L}^2(Y, {\Cal E})$ norm.

									 We thus obtain
									 $$
									 N(\sigma-\sigma_s)^2 \leq \Vert
									 D(\sigma) \Vert^2_{X-B_s}+ \Vert
									 \sigma \Vert ^2_{X-B_s}
									 \frac{M^2}{s^2} \Vert \sigma
									 \Vert^2_X,
									 $$
									 for any $s=1,\dots, \ell$.
									 This completes the proof of Lemma
									 3.3. $\qed$.

									 To end this section, we would like
									 to state separately a very useful
									 fact shown in the proof of Lemma
									 3.3.

							\proclaim{Proposition 3.4 } Let $X$
									 be a non--compact complete orbifold which 
									  is sufficiently regular at
									 infinity, and let $y_0 \in
									 X-\Sigma(X)$ be a fixed 
									 point of $X$. Then there exists a
									 sequence of continuous functions
									 $b_k$, $k\in \Bbb N$, with
									 \roster
									 \item $ b_k: X\to [0,1]$
									 \item $b_k=1$ on $B_k= \{ y\in X|
									 \rho(y) =d(y,y_0)\leq k\}$. 
									 \item The support of $b_k$ is
									 contained in ${\overline B}_{2k}$.
									 \item The function $b_k$ is differentiable almost everywhere
									 and at points of differentiability
									 we have
									 $$
									 \Vert \nabla (b_k) \Vert^2 \leq
									 \frac{M^2}{k^2}, \quad k \in \Bbb N .
									 $$
									 \endroster
									 \endproclaim

									 \vskip 1em
									  \noindent {\bf 4.  The Square of the Dirac Operator.}
									   \vskip 1em

									   As we have seen in Section 3,
									   there is always a unique, closed,
									   self-adjoint extension of a
									   generalized Dirac operator $D_E$
									   on a complete  orbifold
									   $X$ which is sufficiently regular at infinity. This unique extension will
									   still be called $D_E$ and its
									   domain will be denoted by ${\Cal
									   D}(D_E)$. In particular, for any
									   two sections 
									   $\sigma_1, \sigma_2$ of ${\Cal
									   D}(D_E)$, we have
									   $$
									   (D_E\sigma_1, \sigma_2)=
									   (\sigma_1, D_E \sigma_2).
									   $$ 
									   From this we will derive below
									   that if $\sigma \in {\Cal
									   D}(D_E)$, then 
									   $D_E(\sigma)=0$ if and only if
									   $D_E^2(\sigma)=0$. 

						 \proclaim{Theorem 4.1} Let $X$ be a 
non--compact
		  complete 	orbifold which  is sufficiently
									   regular at infinity, and let $E$
									   be a Hermitian orbibundle  (with connection $\nabla^E$) over $X$.
									   Let 			
$ D_E: {\Cal C}^\infty_c (X,
									   \Cal E ) \to {\Cal
C}^\infty_c (X, \Cal E),
									   $
be the generalized Dirac operator
									   on $X$ with coefficients in $E$.
									   Then
									   $D_E(\sigma)=0$ if and only if
									   $D_E^2(\sigma)=0$ for any $\sigma
									   \in {\Cal D}(D_E)$.\endproclaim

									   \noindent {\it Proof.} Set
									   $D=D_E$ and ${\Cal D}={\Cal
									   D}(D_E)$. The non-trivial part of
									   the proof is to show that
									   $D^2\sigma =0$ implies  $D \sigma
									   =0$.
									   Since $D^2$ is elliptic, the
									   equation $D^2\sigma=0$ implies
									   $\sigma \in {\Cal C}^\infty(X,
									   \Cal E)$. In fact, via a
									   partition of unity, we can
									   consider this equation on a chart
									   of a locally finite orbifold
									   atlas. Then, at this 
									   level, we are dealing with a
									   manifold elliptic operator, and
									   therefore standard local theorems
									   on elliptic operators apply, such as the 
									   smoothness of solutions of
									   elliptic systems we need. 
									   Now choose a sequence $\{b_k\}$,
									   $k\in \Bbb N$, as in the proof of
									   Lemma 3.3 and in Proposition 3.4.
									   Then we have, for any $\sigma \in
									   \Cal C (X, \Cal E)$ with $D\sigma =0$,
									   $$
									   (D^2\sigma, b_k^2 \sigma)=
									   (D\sigma, D(b_k^2 \sigma))
									   $$
									   $$
									   =(D\sigma, 2b_k\nabla(b_k)
									   \sigma+b_k^2 D\sigma)=
									   (b_kD\sigma, 2\nabla(b_k)
									   \sigma+b_k D\sigma),
									   $$
									   since
									   $
									   D(f\sigma)= (\nabla f)\sigma +f
									   D(\sigma)$  for any f$\in 
									   \Cal C^\infty (X)$.  almost
									   everywhere.
									   Now we have (recall that $\Vert \,
									   \Vert_X$ is the ${\Cal L}^2(X, \Cal
									  E)$  norm),
									  $$
\Vert  b_k D(\sigma)\Vert^2_X = -2
									   (b_kD\sigma, \nabla(b_k) \sigma)
									   \leq 
									   \frac{M}k (\Vert  b_k
									   D(\sigma)\Vert^2_X + \Vert
									   \sigma\Vert^2_X)
									   $$
									   by Proposition 3.4 and the
									   Schwartz inequality.
									   Since the limit of $\frac{M}k
									   \Vert  \sigma \Vert^2_X$ tends to 0
									   as $k\to +\infty$, it follows
									   that  $\Vert  b_k
									   D(\sigma)\Vert^2_X $ tends to 0 as
									   $k\to +\infty$.
									   But as in the proof of Lemma 3.3,

									   $$
									   \Vert
									   D(\sigma-\sigma_k)\Vert^2_{X}\leq
									   \Vert D(\sigma)\Vert^2_{X-B_k}+
									   \frac{M^2}{k^2} \Vert \sigma
									   \Vert^2_X.
									   $$
									   Hence
									   $$
									   \lim_{k\to +\infty} \Vert
									   D(\sigma-\sigma_k)\Vert^2_{X} =0
									   $$
									   as the union of all $B_k$'s is
									   $X$, and $D(\sigma) \in {\Cal
									   L}^2 (X, \Cal E)$.
									   Then 
									   $$
									   \lim_{k\to +\infty} D(\sigma_k)=
									   D(\sigma)
									   $$
									   in ${\Cal L}^2 (X, \Cal E)$,
									   which implies 
									   $$
									   \lim_{k\to +\infty}\Vert
									   D(\sigma_k) \Vert_X = \Vert
									   D(\sigma) \Vert_X.
									   $$
									   As $D(b_k\sigma)=
									  ( \nabla b_k)\sigma +b_k D(\sigma)$
									   almost everywhere, $\forall k\in \Bbb N$,
									   we have
									   $$
									   \lim_{k\to +\infty}\Vert
									   D(\sigma_k) \Vert_X = \lim_{k\to
									   +\infty}\Vert b_kD(\sigma) \Vert_X
									   $$
									   by Proposition 3.4. Hence
									   $$
									   \Vert D(\sigma) \Vert_X =\lim_{k\to
									   +\infty}\Vert D(\sigma_k) \Vert_X =
									   \lim_{k\to +\infty}\Vert
									   b_kD(\sigma) \Vert_X =0. \quad \qed
									   $$ 
The statement of Theorem 4.1 is
									   also true for sections in ${\Cal
									   L}^2(X, \Cal E)$, as can be shown
									   using approximation.

											 \vskip 1em
											 \noindent {\bf 5.  The Stokes'/
											 Divergence Theorem
											 on Non--Compact Orbifolds.}
											 \vskip 1em

											 In this section we will
											 state and prove a
											 Stokes'/Divergence theorem
											 which is a generalization
											 of manifold results of Gaffney, Yau,
											 and Karp, see [Gn2], [Y],
											 [K]. Our presentation
follows the outline given
											 in [K] for the	 corresponding manifold
											 case. The proof of our
											 theorem relies heavily on
											 the results we proved in
											 Sections 3 and 4.

											 Given a vector field $V$ on
											 an orbifold $X$, choose a
											 locally finite orbifold
											 atlas $\Cal F$,  $\Cal F
											 =\{ ({\tilde U_i}, G_i) | i
											 \in {\Bbb N\} }$, with
											 associated smooth partition
											 of unity $\eta=\{ \eta_i
											 \}_{i\in \Bbb N}$. Then the
											 divergence of $V$ is given
											 in local charts by, [Ch] (here
														  ${\tilde x}=
														  ({\tilde x}_1,
														  \dots, {\tilde
														  x}_n) $
														  denotes the
															  coordinate
															  in
															  ${\tilde
															  U_i}$),
											 $$
											 div(\tilde V)= \sum_k\frac1{\sqrt{\tilde g}} \frac{\partial}{\partial {\tilde x}_k}\left({\sqrt{\tilde g}} {\tilde V}_k\right), \quad {\tilde V}=\sum_k {\tilde V}_k
											 \frac{\partial}{\partial
											 {\tilde x}_k}  \hbox{ on }
											 {\tilde U}_i. 
											 $$

											 \proclaim{Theorem 5.1} Let $X$ be an even--dimensional  non--compact complete $Spin^c$ almost complex orbifold which is sufficiently regular at infinity.
Assume that a connection is chosen on the dual of its canonical line bundle. Let $V$ be a vector field on $X$ such that 
$$
\lim_{k\to +\infty} \inf
\frac1k \int_{B_{2k}-B_k}
\Vert V\Vert\, dv =0,
 $$
where $\Vert V\Vert$
denotes the length of $V$, and $B_k= \{ y\in X|
\rho(y) =d(y,y_0)\leq k\}$ for a fixed $y_0 \in X-\Sigma(X)$, where $\Sigma(X)$ is the singular locus of $X$. 
 Then if either $(div\,
V)^+$ or $(div\, V)^-$ is
integrable on $X$, we have
$$
\int_X div\,(V) \, dv =0.
$$ 
\endproclaim

											  \noindent {\it Proof.}
											  Choose a sequence
											  $\{b_k\}$, $k\in \Bbb N$,
											  as in Proposition 3.4.
											  Integrating $div \, (b_k^2
											  V)$ over $B_{2k}$, for a
											  sufficiently large $k$,
											  and applying the
											  divergence theorem for finite domains, we
											  obtain
											  $$
											  0= \int_ {B_{2k}} div\,
											  (b_k^2\, V) \,dv.
											  $$
											  Hence, by Proposition 3.4,

											  $$
											  \left| \int_{B_{2k}} b_k^2
											  div \,(V)\, dv \right| \leq
											  \frac{M}k
											  \int_{B_{2k}-B_k} \Vert
											  V \, \Vert\, dv.
											  $$
											  Thus, if we for example
											  suppose  $(div \, V)^-$
											  integrable, the above
											  inequality implies
											  $$
											  \int_{B_{k}} (div \,(V))^+
											  dv -\int_X (div \,(V))^-
											  dv \leq \frac{M}k
											  \int_{B_{2k}-B_k} \Vert
											  V \Vert \, dv.
											  $$
											  Because of our hypothesis,
											  we can choose a sequence
											  $k(j)
											  \to + \infty$, such that 
											  $$
											  \lim_{j\to + \infty
											  }\int_{B_{2k(j)}-B_k(j)} \Vert
											  V\Vert\, dv  =0.
											  $$
											  Consequently, $(div
											  \,(V))^- $ is also
											  integrable, and 
											  $$
											  \int_X div \,(V) dv \leq
											  0.
											  $$
											  But now the same argument
											  can be repeated started
											  from  $(div \,(V))^- $.
											  Hence
											  $$
	\int_X div \,(V) dv = 0.  \quad \qed
											  $$

											  \proclaim{Corollary 5.2}
											  Let $X$ and $V$ be as in
											  Theorem 5.1, and also
											  assume that $X$ has q-th
											  order volume growth (i.e.,
											  there exists $c>0$ and
											  $q\geq 1$ such that $vol\,
											  (B_k) \leq c \,k^q$,
											  $\forall k\geq 1$.) If
											  $div\, (V)\geq 0$ outside
											  of some compact set, and
											  either 
											  \roster
											  \item $q>1$ and $V\in
											  {\Cal L}^p(X, \Cal E)$,
											  with $\frac1p + \frac 1q
											  =1$, or
											  \item $q=1$ and $\Vert
											  V\Vert \to 0$ uniformely
											  at $\infty$ in $X$, then 
											  \endroster  

											  $$
											  \int_X div \,(V) dv = 0.
											  $$
											  \endproclaim

											  \noindent {\it Proof.}
											  Very similar to [K; Proof
											  of Corollary 1]. $\qed$

											  \vskip 1em
											  \noindent {\bf 6.  Some
											  Vanishing Theorems.}
											  \vskip 1em

											  The results in this
											  section are a
											  generalization to
											  orbifolds of some of the
											  results proved for
											  manifolds by Gromov
											  and Lawson in [GL]. In this section, we will
											  let $ E =\Bbb C$
											  unless otherwise noticed.
											  We will also substitute
											  $\Cal S$ for $\Cal E$. 

											  The scalar orbifold Laplacian $\Delta$ can
											  be defined in analogy with
											  the Laplacian on
											  manifolds. (c.f. [Ch]
											  Section 2.) In fact, on an
											  orbifold chart ${\tilde U}_i$
											  of a
											  standard orbifold atlas
											  $\Cal U  =\{ ({\tilde U_i},
											  G_i )| i \in { I\} }$ of
											  $X$, we define, (here
											  ${\tilde x}= ({\tilde x}_1, \dots, {\tilde x}_n) $ denotes the
											  coordinate in  ${\tilde
											  U_i}$),
											  $$
											  \Delta \tilde u = \sum_{k,j}
											  -\tilde{g}^{k,j}
											  \frac{\partial^2{\tilde
											  u}}{\partial{\tilde x}_k
											  \partial{\tilde x}_j}
											  -\sum_j{\tilde B}_j
											  \frac{\partial {\tilde
											  u}}{\partial{\tilde x}_j},
											  $$  
											  with
											  $$
											  {\tilde B}_j=
											  \frac1{2\tilde{g}}
											  \sum_k
											  \frac{\partial{\tilde{g}}}{\partial{\tilde x}_k}{\tilde{g}^{k,j}}+
											  \sum_k
											  \frac{\partial{\tilde{g_{k,j}}}}{\partial{\tilde x}_k}, \quad \forall
											  u \in {\Cal
											  C}^\infty(U_i).
											  $$
											  In the above expression,
											  $\tilde{g}={\tilde{g}^{k,j}}$,
											  $k,j=1,\dots, n$, is a
											  $G_i$--invariant metric.
											  Laplacians can also be
											  defined to act on general orbibundles such as $\Cal
											  S$ by using the above
											  definition on										  orbisections. The
											  following Green's formula
											  holds.

											  \proclaim{Proposition 6.1}
											  Let $X$ be a non--compact
											  complete orbifold which is 
											  sufficiently regular at
											  infinity, and let $\Cal S$ be  the $Spin^c$ bundle of $X$. 
											  Then for any two
											  sections $\sigma_j$,
											  $j=1,2$ in ${\Cal
											  C}^\infty (X, \Cal S)$, at
											  least one of which with
											  compact support, we have
											  $$
											  \int_X <\Delta \sigma_1,
											  \sigma_2 > dv =\int_X <
											  \nabla\sigma_1,
											  \nabla\sigma_2 > dv 
											  $$ 
											  \endproclaim

											  \noindent {\it Proof.}
											  Because of our hypothesis
											  at infinity, the proof given in [Ch; Section 2] in the scalar case is also valid here. In particular, this result
											  is an orbifold version of
											  [Si; Proposition 1.2.2],
											  which can be proved as in
											  the manifold case.
											  $\qed$

											  Proposition 6.1 motivates
											  us to choose the Sobolev
											  norm
											  $$
											  \Vert \sigma\Vert_1^2=
											  \int_X\left( <\sigma,
											  \sigma>+ <\nabla \sigma,
											  \nabla \sigma>\right)dv =
											  \int_X\left( <\sigma,
											  \sigma>+ <\Delta\sigma,
											  \sigma>
											  \right)dv
											  $$
											  Thus, by reasoning as in
											  Section 5, we obtain,

											  \proclaim{Theorem 6.2} Let $X$ be a non--compact complete 
 orbifold which  is sufficiently regular at infinity. Then the domain
 of the unique closed
 self-adjoint extension of the $Spin^c$ Laplacian
 $\Delta$,  $  \Delta: {\Cal C}^\infty_c
 (X, {\Cal S}) \to {\Cal
 L}^2 (X, {\Cal S}), $
 is the completion $ {\Cal L}^{1,2}(X, {\Cal
 S})$ of ${\Cal C}^\infty_c (X, {\Cal S}) $ in the norm 
 $\Vert \,\Vert_1$.
 Furthermore,
$\Delta(\sigma)=0$ if and
											  only if										  $\nabla(\sigma)=0$, i.e.,
											  $\sigma$ is parallel.  
											  \endproclaim

											  \noindent {\it Proof.} We
											  only need to prove the
											  last claim, which follows
											  from Proposition 6.1 in
											  the case of sections with
											  compact support. The
											  general case follows from
											  Laplacian analogs of
											  Lemmas 3.2 and 3.3. 
											  $\qed$

											  The important
											  Bochner--Weitzenbr\"ock
											  formula, a classic result
											  for manifolds, can also be
											  easily extended to
											  orbifolds, by using local
											  coordinates. 

											  \proclaim{Proposition 6.3}  Let $X$ be a   non--compact
complete  orbifold which  is	 sufficiently regular at  infinity.  If $D$ is the Dirac operator on $X$ with
 coefficients in the $Spin^c$  bundle $\Cal S$, and $\Delta$ is the $Spin^c$ Laplacian, then
 $$
D^2= \Delta+ {\Cal R},
$$
where ${\Cal R}$ is given below (c.f. [Du; Theorem 6.1], [LM; Theorem D12] for the manifold case), 
$$
\Cal R= \frac14\,  k\,  + \frac12 c(K^*),
$$
where $k$ is the scalar curvature,  and $ c(K^*)$ denotes the
Clifford multiplication of the curvature 2 form of the fixed  
connection on the line bundle $K^*$.
\endproclaim

As a consequence of the
above formula, we obtain,
as in [GL; Theorem 2.8],

											  \proclaim{Theorem 6.4} Let
											  $X$ be a non--—compact complete
											  orbifold which  is
											  sufficiently regular at
											  infinity. If $D$ is the
											  Dirac operator on $X$ with
											  coefficients in the
											  $Spin^c$ bundle $\Cal S$,
											  then
											  the domain $\Cal D$ of the
											  unique self-adjoint
											  extension of $D$ is
											  exactly
											  $$
											  {\Cal L}^{1,2}(X, {\Cal
											  S}), \quad \hbox{that is,}
											  $$
											  the completion of ${\Cal
											  C}^\infty_c (X, {\Cal S})
											  $  in the norm 
											  $$
											  \Vert \sigma\Vert_1^2=
											  \int_X\left( <\sigma,
											  \sigma>+ <\nabla\sigma,
											  \nabla\sigma>\right)dv =
											  \int_X\left( <\sigma,
											  \sigma>+ <\Delta_\sigma,
											  \sigma>
											  \right)dv
											  $$
											  Moreover, for every
											  $\sigma\in \Cal D$,
											  $$
											  \Vert D\sigma\Vert^2_X =
											  \Vert \nabla\sigma \Vert^2_X
											  + ({\Cal R} \sigma,
											  \sigma), 
											  $$
 where $\Vert\, \Vert_X$  denotes the ${\Cal L}^2$	 norm, 
$\Cal R$ is as in Theorem 6.3, 
and $(\, , )$ th ${\Cal L}^2$ inner  product. \endproclaim
											  \noindent {\it Proof.} For
 sections with compact
											  support Theorem 6.4
											  follows directly from the
											  Bochner--Weitzenbr\"ock
											  formula and the
											  self-adjointness of the
											  Dirac operator. More in
  general, approximate a
  section $\sigma \in {\Cal
											  D}$, via the ${\Cal L}^2$
											  norm, by a sequence
											  $\sigma_k$ of compact
											  support
											  such that $D\sigma_k \to
											  D\sigma$.  (This is
											  possible because $\Cal D$
											  is in particular equal to
											  the minimal domain of
											  $D$.) By passing to the
											  limit, we obtain, 
											  $$
											  \Vert D\sigma\Vert^2_X =
											  \Vert\nabla\sigma \Vert^2_X
											  +({\Cal R}\sigma, \sigma),
											  $$
											  since   ${\Cal R}$ is
											  bounded. $\qed$

											  The following corollaries
											  can be derived as in the
											  manifold case
											  (see [GL; Section 2]). We
	 will thus leave their
	  proofs to the reader.

											  \proclaim{Corollary 6.5}
Let $X$ be a 
non--compact, complete
  orbifold which  is
											  sufficiently regular at
											  infinity. Let $D$ and
											  $\Cal R$ be as in 
											  Theorem 6.4. Suppose that
											  $\Cal R > 0 $ pointwise on
											  $X$. Then
											  $$
											  \hbox {Ker} (D) =\hbox
											  {Coker} (D) =0.
											  $$
											  If furthermore, $\Cal
											  R\geq\, c \, Id$, for some
											  constant $c>0$, then 
											  $D: {\Cal L}^{1,2}(X,
											  {\Cal S}) \to {\Cal
											  L}^{2}(X, {\Cal S})$  
											  is an isomorphism of
											  Hilbert spaces. In this
											  case, $D^{-1}: {\Cal
											  L}^{1,2}(X, {\Cal S}) \to
											  {\Cal L}^{2}(X, {\Cal S})$
											  is also a bounded
											  operator.
											  \endproclaim

											  \proclaim{Corollary 6.6}
											  Let $X$ ,  $D$ and  $\Cal
											  R$ be as in 
											  Corollary 6.5, with $\Cal
											  R\geq \,c \,Id$, for some
											  constant $c>0$. Since
$X$ is  even-dimensional,
											  $D: {\Cal L}^{1,2}(X,
											  {\Cal S}) \to {\Cal
											  L}^{2}(X, {\Cal S})  $
											  splits into its $\pm$
											  decomposition, $D^+$ and
											  $D^-$, see Section 1.
Then both $D^+$ and $D^-$
have bounded inverses. \endproclaim

											  \proclaim{Remark 6.7} The
											  results of this section
											  can also be proved, with suitable modifications for 
											  generalized Dirac
											  operators with
											  coefficients in any 
Hermitian 	orbibundle  $E$ (with connection $\nabla^E$).
											  \endproclaim

											  \vskip 1em
											  \noindent {\bf References}
\itemitem {[At1]} M.F. Atiyah, Elliptic operators, discrete groups and von
 Neumann algebras, Colloque "Analyse et Topologie" en	 l'honneur de Henri Cartan (Orsay, 1974),  pp. 43--72. Asterisque, No. 32-33, Soc. Math. France, Paris, 1976.
\itemitem {[AX]} 
B. Apanasov and X.  Xie, 
Discrete actions on nilpotent Lie groups and negatively curved spaces,
Differential Geom. Appl. 20 (2004), 11--29. 
											 \itemitem {[B]} J. Borzellino,  Orbifolds of maximal
diameter, Indiana Univ. Math. J.  42 (1993), 37--53.
			
													   \itemitem {[BGV]} N. Berline, E. Getzler and M. Vergne, Heat kernels and Dirac operators, Grundleheren der Mathematical Wissenshaften 298, Springer--Verlag, Berlin, 1992.
															   \itemitem
															   {[C]} G.
															   Chen,
															   Calculus
															   on
															   orbifolds,
															   Sichuan
															   Daxue
															   Xuebao
															   41
															   (2004),
															   931--939.

															   \itemitem
															   {[Ch]}
															   Y.-C.
															   Chiang,
															   Harmonic
															   maps of
															   $V$-manifolds,           Ann. Global Anal. Geom.  8  (1990), 315-344.
															   \itemitem
															   {[Du]} J.
															   J.
															   Duistermaat,
															   The heat
															   kernel
															   Lefschetz
															   fixed
															   point
															   formula
															   for the
															   $Spin^c$
															   Dirac
															   operator,
															   Progress
															   in
															   Nonlinear
															   Differential
															   Equations
															   and their
															   Applications,
															   18. Birkh\"auser, Inc., Boston,
				 MA, 1996.

\itemitem {[Fa1]} C. Farsi, $K$-theoretical index theorems for orbifolds, Quat.  J.  Math.  43 (92), 183--200.  
\itemitem {[Fa2]} $\underline {\hskip 1.0in}$ Orbifold spectral theory, Rocky Mtn.  J. Math. 31(2001), 215--235.
\itemitem {[Fa3]} $\underline {\hskip 1.0in}$ Orbifold $\eta$-invariants, Indiana Math. Journal, 
Indiana Math. J. 35 (2007),  501-521.
\itemitem {[Fa4]} $\underline {\hskip 1.0in}$ A relative orbifold index theorem,  J. Geom Phys. 8 (2007), 1653-1668.
\itemitem {[GL]} M.  Gromov and M.  Lawson, Positive scalar curvature and the Dirac operator on complete Riemannian manifolds, Inst.  Hautes Etudes Sci.  Publ.  Math.  No.  58, (1983), 83--196.
\itemitem {[Gn1]} M. Gaffney,  The harmonic operator for exterior differential forms, Proc.  Nat.  Acad.  Sci.  U. S.  A.  37 (1951),  48--50. 
\itemitem {[Gn2]} M.  Gaffney, A special Stokes's theorem for complete Riemannian manifolds, Ann.  of Math.  (2) 60 (1954), 140--145. 
\itemitem {[K]} L.  Karp, On Stokes' theorem for noncompact manifolds, Proc.  Amer.  Math.  Soc.  82 (1981), 487--490.  
\itemitem {[Kw1]} T.  Kawasaki, The signature theorem for $V$--manifolds, Topology 17 (78), 75--83.  \itemitem {[Kw2]} 
$\underline {\hskip 1.0in}$ The Riemann Roch theorem for complex $V$--manifolds, Osaka J. Math. 16 (1979), 151--159.
\itemitem {[Kw3]} $\underline {\hskip 1.0in}$ The index of elliptic operators over $V$-manifolds, Nagoya Math. J. 84 (81), 135--157.
\itemitem {[LM]} H. B. Lawson, Jr, and M.-L. Michelsohn, Spin Geometry, Princeton University Press, Princeton, New Jersey, 1989.
\itemitem {[LR]} W.  L\"uck and J.  Rosenberg, Equivariant Euler characteristics and $K$-homology Euler classes for proper cocompact $G$-manifolds, Geom.  Topol.  7 (2003), 569--613. 
\itemitem {[LoR]} D. D. Long, and A. W. Reid,
All flat manifolds are cusps of hyperbolic orbifolds,
Algebr. Geom. Topol. 2 (2002), 285--296. 
\itemitem {[N]} Y. Nakagawa, An isoperimetric inequality for orbifolds, Osaka J. Math. 30 (1993), 733--739. 
\itemitem {[Si]} J.  Simons, Minimal varieties in Riemannian manifolds, Ann.  of Math.  (2) 88 (1968), 62--105.  
\itemitem {[Stan]} E.  Stanhope, Spectral bounds on orbifold isotropy,  Ann.  Global Anal.  Geom.  27 (2005), 355--375. 
\itemitem {[V]} M. Vergne, 
Equivariant index formulas for orbifolds.
Duke Math. J. 82 (1996), no. 3, 637--652.
\itemitem {[W]} J.  Wolf, Essential self-adjointness for the Dirac operator and its square,  Indiana Univ. Math. J.  22 (1972/73), 611--640.
\itemitem {[Y]} S.  T. Yau, Some function-theoretic properties of complete Riemannian manifold and their applications to geometry, Indiana Univ.  Math.  J.  25 (1976), 659--670.
\end